# Asymptotic normality and consistency of a two-stage generalized least squares estimator in the growth curve model

JIANHUA HU[1] and GUOHUA YAN[2]

[1]*Department of Statistics, Shanghai University of Finance and Economics, 777 Guo Ding Road, Shanghai 200433, P.R. China. E-mail: frank.jianhuahu@gmail.com*
[2]*Department of Statistics, University of British Columbia, 333-6356 Agricultural Road, Vancouver, BC, Canada, V6T 1Z2. E-mail: guohuayan@gmail.com*

Let $\mathbf{Y} = \mathbf{X\Theta Z'} + \boldsymbol{\mathcal{E}}$ be the growth curve model with $\boldsymbol{\mathcal{E}}$ distributed with mean $\mathbf{0}$ and covariance $\mathbf{I}_n \otimes \boldsymbol{\Sigma}$, where $\boldsymbol{\Theta}$, $\boldsymbol{\Sigma}$ are unknown matrices of parameters and $\mathbf{X}$, $\mathbf{Z}$ are known matrices. For the estimable parametric transformation of the form $\boldsymbol{\gamma} = \mathbf{C\Theta D'}$ with given $\mathbf{C}$ and $\mathbf{D}$, the two-stage generalized least-squares estimator $\hat{\boldsymbol{\gamma}}(\mathbf{Y})$ defined in (7) converges in probability to $\boldsymbol{\gamma}$ as the sample size $n$ tends to infinity and, further, $\sqrt{n}[\hat{\boldsymbol{\gamma}}(\mathbf{Y}) - \boldsymbol{\gamma}]$ converges in distribution to the multivariate normal distribution $\mathcal{N}(\mathbf{0}, (\mathbf{CR}^{-1}\mathbf{C'}) \otimes (\mathbf{D}(\mathbf{Z'\Sigma}^{-1}\mathbf{Z})^{-1}\mathbf{D'}))$ under the condition that $\lim_{n \to \infty} \mathbf{X'X}/n = \mathbf{R}$ for some positive definite matrix $\mathbf{R}$. Moreover, the unbiased and invariant quadratic estimator $\hat{\boldsymbol{\Sigma}}(\mathbf{Y})$ defined in (6) is also proved to be consistent with the second-order parameter matrix $\boldsymbol{\Sigma}$.

## 1. Introduction

The growth curve model is defined as

$$\mathbf{Y} = \mathbf{X\Theta Z'} + \boldsymbol{\mathcal{E}}, \qquad \boldsymbol{\mathcal{E}} \sim \mathcal{G}(\mathbf{0}, \mathbf{I}_n \otimes \boldsymbol{\Sigma}), \tag{1}$$

where $\mathbf{Y}$ is an $n \times p$ matrix of observations, $\mathbf{X}$ and $\mathbf{Z}$ are known $n \times m$ ($n > m$) and $p \times q$ ($p > q$) full-rank design matrices, respectively, $\boldsymbol{\Theta}$ is an unknown $m \times q$ matrix, called the first-order parameter matrix, and $\boldsymbol{\Sigma}$ is an unknown positive definite matrix of order $p$, called the second-order parameter matrix. $\boldsymbol{\mathcal{E}}$ follows a general continuous distribution $\mathcal{G}$ with mean matrix $\mathbf{0}$ and Kronecker product structure covariance matrix $\mathbf{I}_n \otimes \boldsymbol{\Sigma}$.

Model (1) was proposed by Potthoff and Roy [11] under the normality assumption of the error matrix $\boldsymbol{\mathcal{E}}$. Since then, parameter estimation, hypothesis testing and prediction of future values have been investigated by numerous researchers, generating a substantial amount of literature concerning the model.

In what follows, we give a brief review of the literature on large sample properties for the growth curve model, a particular kind of multivariate regression model. Chakravorti [2] presented the asymptotic properties of the maximum likelihood estimators. Žežula [15] investigated the asymptotic properties of the growth curve model with covariance components. Gong [4] gave the asymptotic distribution of the likelihood ratio statistic for testing sphericity. Bischoff [1] considered some asymptotic optimal tests for some growth curve models under non-normal error structure. However, no work has been done on the asymptotic normality and consistency of two-stage generalized least-squares estimators of the first-order parameter matrix for the growth curve model (1).

In this paper, we shall investigate the consistency and asymptotic normality of a two-stage generalized least-squares estimator $\hat{\boldsymbol{\gamma}}(\mathbf{Y})$ for the estimable parametric transformation of the form $\boldsymbol{\gamma} = \mathbf{C}\boldsymbol{\Theta}\mathbf{D}'$ with respect to the first-order parameter matrix $\boldsymbol{\Theta}$. In addition, we shall demonstrate the consistency of a known quadratic covariance estimator $\hat{\boldsymbol{\Sigma}}(\mathbf{Y})$ with the second-order parameter matrix $\boldsymbol{\Sigma}$ (see Žežula [14]).

Readers are referred to Eicker [3], Theil [13] and Nussbaum [10] for results on the large sample properties of the least-squares estimators for ordinary univariate and multivariate regression models.

This paper is divided into four sections. Some preliminaries are presented in Section 2. In particular, for the estimable parametric transformation of the form $\boldsymbol{\gamma} = \mathbf{C}\boldsymbol{\Theta}\mathbf{D}'$, a two-stage generalized least-squares estimator $\hat{\boldsymbol{\gamma}}(\mathbf{Y})$ is defined in (7). The consistency of the estimator $\hat{\boldsymbol{\gamma}}(\mathbf{Y})$ and the consistency of the known quadratic estimator $\hat{\boldsymbol{\Sigma}}(\mathbf{Y})$ defined in (6) are investigated in Section 3. Finally, in Section 4, the asymptotic normality of the two-stage generalized least-squares estimator $\hat{\boldsymbol{\gamma}}(\mathbf{Y})$ is obtained under a certain condition.

## 2. Preliminaries

Throughout this paper, the following notation is used. Let $\mathscr{M}_{n \times p}$ denote the set of all $n \times p$ matrices. Let $\mathbf{A}'$ denote the transpose of the matrix $\mathbf{A}$. Let $\operatorname{tr}(\mathbf{A})$ denote the trace of the matrix $\mathbf{A}$. Let $\mathbf{I}_n$ denote the identity matrix of order $n$. For a sequence of numbers $\{a_n\}$ and a sequence of numbers $\{b_n\}$, we say that $a_n = O(b_n)$ if there is a constant $c$ such that $\limsup |a_n/b_n| \leq c$; we say that $a_n = o(b_n)$ if $\lim a_n/b_n = 0$. For an $n \times p$ matrix $\mathbf{Y}$, we write $\mathbf{Y} = [\mathbf{y}_1, \mathbf{y}_2, \ldots, \mathbf{y}_n]'$, $\mathbf{y}_i \in \Re^p$, where $\Re^p$ is the $p$-dimensional real space, and $\operatorname{vec}(\mathbf{Y}')$ denotes the $np$-dimensional vector $[\mathbf{y}_1', \mathbf{y}_2', \ldots, \mathbf{y}_n']'$. Here, the vec operator transforms a matrix into a vector by stacking the columns of the matrix one under another. $\mathbf{Y} \sim \mathcal{G}(\mathbf{M}, \mathbf{I}_n \otimes \boldsymbol{\Sigma})$ means that $\mathbf{Y}$ follows a general continuous distribution $\mathcal{G}$ with $\mathrm{E}(\mathbf{Y}) = \mathbf{M}$ and that $\mathbf{I}_n \otimes \boldsymbol{\Sigma}$ is the covariance matrix of the vector $\operatorname{vec}(\mathbf{Y}')$; see Muirhead [9], Section 3.1. The Kronecker product $\mathbf{A} \otimes \mathbf{B}$ of matrices $\mathbf{A}$ and $\mathbf{B}$ is defined to be $\mathbf{A} \otimes \mathbf{B} = (a_{ij}\mathbf{B})$. We then have $\operatorname{vec}(\mathbf{ABC}) = (\mathbf{C}' \otimes \mathbf{A})\operatorname{vec}(\mathbf{B})$. Let $\mathbf{A}^+$ denote the Moore–Penrose inverse of $\mathbf{A}$ and $\mathbf{P}_X = \mathbf{X}(\mathbf{X}'\mathbf{X})^{-}\mathbf{X}'$ be the projection onto the column space $\mathscr{C}(\mathbf{X})$ of a matrix $\mathbf{X}$ along the orthogonal complement $\mathscr{C}(\mathbf{X})^{\perp}$ of $\mathscr{C}(\mathbf{X})$.

Given $\mathbf{A} \in \mathcal{M}_{n \times p}$ and $\mathbf{B} \in \mathcal{M}_{p \times s}$, a linear parametric function $\mathbf{B}'\boldsymbol{\beta}$ is called *estimable* with respect to $\mathbf{A}$ if there exists some $\mathbf{T} \in \mathcal{M}_{n \times s}$ such that $\mathrm{E}(\mathbf{T}'\mathbf{A}\boldsymbol{\beta}) = \mathbf{B}'\boldsymbol{\beta}$ for all $\boldsymbol{\beta} \in \Re^p$; see Hu and Shi [6] for a more detailed description.

Note that the first-order parameter $\boldsymbol{\Theta}$ in model (1) is defined before a design is planned and observation $\mathbf{Y}$ is obtained. Thus, the rows of the design matrix $\mathbf{X}$ in model (1) are added one after another and the term $\mathbf{Z}$ in model (1) does not depend on the sample size $n$; see the example in Potthoff and Roy [11]. So, we shall only consider the case of full-rank matrices $\mathbf{X}$ and $\mathbf{Z}$.

As discussed in Potthoff and Roy [11], hypotheses of the form $\mathbf{C}\boldsymbol{\Theta}\mathbf{D}' = \mathbf{0}$ under model (1) are usually considered, where $\mathbf{C} \in \mathcal{M}_{s \times m}$ and $\mathbf{D} \in \mathcal{M}_{t \times q}$. Thus, in this paper, we shall consider the estimator of the parametric transformation $\boldsymbol{\gamma} = \mathbf{C}\boldsymbol{\Theta}\mathbf{D}'$ of $\boldsymbol{\Theta}$ with given matrices $\mathbf{C} \in \mathcal{M}_{s \times m}$ and $\mathbf{D} \in \mathcal{M}_{t \times q}$.

We shall begin by reviewing the case of a known second-order parameter matrix $\boldsymbol{\Sigma}$, say $\boldsymbol{\Sigma}_0$. According to the theory of least squares (see, e.g., Rao [12], 4a.2), the normal equations of model (1) are $\mathbf{X}'\mathbf{X}\boldsymbol{\Theta}\mathbf{Z}'\boldsymbol{\Sigma}_0^{-1}\mathbf{Z} = \mathbf{X}'\mathbf{Y}\boldsymbol{\Sigma}_0^{-1}\mathbf{Z}$. The least-squares estimator $\hat{\boldsymbol{\Theta}}_0$ of $\boldsymbol{\Theta}$ is given by

$$\hat{\boldsymbol{\Theta}}_0 = (\mathbf{X}'\mathbf{X})^{-1}\mathbf{X}'\mathbf{Y}\boldsymbol{\Sigma}_0^{-1}\mathbf{Z}(\mathbf{Z}'\boldsymbol{\Sigma}_0^{-1}\mathbf{Z})^{-1}. \tag{2}$$

Since

$$(\mathbf{X}'\mathbf{X})^{-1}\mathbf{X}' = (\mathbf{X}'\mathbf{X})^{-1}\mathbf{X}'\mathbf{P}_X$$

and

$$\mathbf{Z}(\mathbf{Z}'\boldsymbol{\Sigma}_0^{-1}\mathbf{Z})^{-1}\mathbf{Z}' = (\mathbf{P}_Z\boldsymbol{\Sigma}_0^{-1}\mathbf{P}_Z)^{+}, \tag{3}$$

(2) can be written as

$$\hat{\boldsymbol{\Theta}}_0 = (\mathbf{X}'\mathbf{X})^{-1}\mathbf{X}'\mathbf{P}_X\mathbf{Y}\boldsymbol{\Sigma}_0^{-1}(\mathbf{P}_Z\boldsymbol{\Sigma}_0^{-1}\mathbf{P}_Z)^{+}\mathbf{Z}(\mathbf{Z}'\mathbf{Z})^{-1}. \tag{4}$$

Let

$$\hat{\boldsymbol{\gamma}}_0 = \mathbf{C}\hat{\boldsymbol{\Theta}}_0\mathbf{D}'. \tag{5}$$

The mean and covariance of $\hat{\boldsymbol{\gamma}}_0$ are, respectively, $\mathbf{C}\boldsymbol{\Theta}\mathbf{D}'$ and $(\mathbf{C}(\mathbf{X}'\mathbf{X})^{-1}\mathbf{C}') \otimes (\mathbf{D}(\mathbf{Z}'\boldsymbol{\Sigma}_0^{-1}\mathbf{Z})^{-1}\mathbf{D}')$.

In addition, it follows from Rao [12], 4a.2, that $\boldsymbol{\gamma} = \mathbf{C}\boldsymbol{\Theta}\mathbf{D}'$, for any matrices $\mathbf{C} \in \mathcal{M}_{s \times m}$ and $\mathbf{D} \in \mathcal{M}_{t \times q}$, is an estimable parametric transformation if matrices $\mathbf{X}$ and $\mathbf{Z}$ are of full rank. So, $\hat{\boldsymbol{\gamma}}_0$ defined in (5) is said to be a least-squares estimator of the estimable parametric transformation $\boldsymbol{\gamma} = \mathbf{C}\boldsymbol{\Theta}\mathbf{D}'$. It is easily derived from 4a.2 of Rao [12] that $\hat{\boldsymbol{\gamma}}_0$ is the best linear unbiased estimator (BLUE) of $\boldsymbol{\gamma}$.

Now, we shall focus our attention on the case of an unknown $\boldsymbol{\Sigma}$.

Let

$$\hat{\boldsymbol{\Sigma}}(\mathbf{Y}) = \mathbf{Y}'\mathbf{W}\mathbf{Y}, \qquad \mathbf{W} \equiv \frac{1}{n - \mathrm{rank}(\mathbf{X})}(\mathbf{I} - \mathbf{P}_X). \tag{6}$$

It is well known that $\hat{\boldsymbol{\Sigma}}(\mathbf{Y})^{-1}$ is positive definite with probability 1 (see the proof of Muirhead [9], Theorem 3.1.4). Žežula [14], Theorem 3.7, tells us that $\hat{\boldsymbol{\Sigma}}(\mathbf{Y})$ is a uniformly minimum variance unbiased invariant estimator of $\boldsymbol{\Sigma}$ under the assumption of normality. This estimator $\hat{\boldsymbol{\Sigma}}(\mathbf{Y})$ is often used to find the first-stage estimator; see, for example, Žežula [16]. We shall also take the estimator as the first-stage estimator in our following discussion.

In (5), an unbiased least-squares estimator of $\boldsymbol{\gamma}$ is given when $\boldsymbol{\Sigma}$ is known. However, when $\boldsymbol{\Sigma}$ is unknown, if we write $\hat{\boldsymbol{\Theta}} \equiv (\mathbf{X}'\mathbf{X})^{-1}\mathbf{X}'\mathbf{Y}\boldsymbol{\Sigma}^{-1}\mathbf{Z}(\mathbf{Z}'\boldsymbol{\Sigma}^{-1}\mathbf{Z})^{-1}$, then the statistic $\hat{\boldsymbol{\gamma}} \equiv \mathbf{C}\hat{\boldsymbol{\Theta}}\mathbf{D}'$ depends on $\boldsymbol{\Sigma}$. In this case, we shall use a method called two-stage estimation to find an estimator, which is denoted by $\hat{\boldsymbol{\gamma}}(\mathbf{Y})$: first, based on data $\mathbf{Y}$, we find a first-stage estimator $\tilde{\boldsymbol{\Sigma}}$ of $\boldsymbol{\Sigma}$; second, replace the unknown $\boldsymbol{\Sigma}$ with the first-stage estimator $\tilde{\boldsymbol{\Sigma}}$ and then find $\hat{\boldsymbol{\Theta}}$ through the normal equations of model (1).

We take $\hat{\boldsymbol{\Sigma}}(\mathbf{Y})$ in (6) as the first-stage estimator $\tilde{\boldsymbol{\Sigma}}$. Replacing $\boldsymbol{\Sigma}$ in (4) with $\hat{\boldsymbol{\Sigma}}(\mathbf{Y})$, (5) can be expressed as

$$\hat{\boldsymbol{\gamma}}(\mathbf{Y}) = \mathbf{C}(\mathbf{X}'\mathbf{X})^{-1}\mathbf{X}'\mathbf{Y}\hat{\boldsymbol{\Sigma}}^{-1}(\mathbf{Y})\mathbf{Z}(\mathbf{Z}'\hat{\boldsymbol{\Sigma}}^{-1}(\mathbf{Y})\mathbf{Z})^{-1}\mathbf{D}'. \tag{7}$$

Let

$$\mathbf{H}(\mathbf{Y}) \equiv \hat{\boldsymbol{\Sigma}}^{-1}(\mathbf{Y})(\mathbf{P}_Z\hat{\boldsymbol{\Sigma}}^{-1}(\mathbf{Y})\mathbf{P}_Z)^+. \tag{8}$$

Then, by (3), (7) can be rewritten as

$$\hat{\boldsymbol{\gamma}}(\mathbf{Y}) = \mathbf{C}(\mathbf{X}'\mathbf{X})^{-1}\mathbf{X}'\mathbf{Y}\mathbf{H}(\mathbf{Y})\mathbf{Z}(\mathbf{Z}'\mathbf{Z})^{-1}\mathbf{D}'. \tag{9}$$

The estimator $\hat{\boldsymbol{\gamma}}(\mathbf{Y})$ is said to be a two-stage generalized least-squares estimator of the estimable parametric transformation $\boldsymbol{\gamma} = \mathbf{C}\boldsymbol{\Theta}\mathbf{D}'$.

In the special case of $\mathbf{C}$ and $\mathbf{D}$ being identity matrices, the estimable parametric transformation $\boldsymbol{\gamma}$ is the first-order parameter matrix $\boldsymbol{\Theta}$. By (9) or (4), we have

$$\hat{\boldsymbol{\Theta}}(\mathbf{Y}) = (\mathbf{X}'\mathbf{X})^{-1}\mathbf{X}'\mathbf{Y}\mathbf{H}(\mathbf{Y})\mathbf{Z}(\mathbf{Z}'\mathbf{Z})^{-1}. \tag{10}$$

The following lemma concerns the unbiasedness of the estimator $\hat{\boldsymbol{\gamma}}(\mathbf{Y})$ under the assumption that $\boldsymbol{\mathcal{E}}$ is symmetric about the origin.

**Lemma 2.1.** *Assume that the distribution of $\boldsymbol{\mathcal{E}}$ is symmetric about the origin. Then the statistic $\hat{\boldsymbol{\gamma}}(\mathbf{Y})$ defined in (9) is an unbiased estimator of the estimable parametric transformation $\boldsymbol{\gamma}$.*

**Proof.** Since $\hat{\boldsymbol{\Sigma}}(\mathbf{Y}) = \hat{\boldsymbol{\Sigma}}(\boldsymbol{\mathcal{E}}) = \hat{\boldsymbol{\Sigma}}(-\boldsymbol{\mathcal{E}})$, $\hat{\boldsymbol{\gamma}}(\mathbf{Y})$ can be expressed as

$$\hat{\boldsymbol{\gamma}}(\mathbf{Y}) = \mathbf{C}(\mathbf{X}'\mathbf{X})^{-1}\mathbf{X}'\mathbf{X}\boldsymbol{\Theta}\mathbf{Z}'\hat{\boldsymbol{\Sigma}}^{-1}(\boldsymbol{\mathcal{E}})\mathbf{Z}(\mathbf{Z}'\hat{\boldsymbol{\Sigma}}^{-1}(\boldsymbol{\mathcal{E}})\mathbf{Z})^{-1}\mathbf{D}'$$
$$+ \mathbf{C}(\mathbf{X}'\mathbf{X})^{-1}\mathbf{X}'\boldsymbol{\mathcal{E}}\hat{\boldsymbol{\Sigma}}^{-1}(\boldsymbol{\mathcal{E}})\mathbf{Z}(\mathbf{Z}'\hat{\boldsymbol{\Sigma}}^{-1}(\boldsymbol{\mathcal{E}})\mathbf{Z})^{-1}\mathbf{D}'.$$

Let

$$\mathbf{M}(\boldsymbol{\mathcal{E}}) = \mathbf{C}(\mathbf{X}'\mathbf{X})^{-1}\mathbf{X}'\boldsymbol{\mathcal{E}}\hat{\boldsymbol{\Sigma}}^{-1}(\boldsymbol{\mathcal{E}})\mathbf{Z}(\mathbf{Z}'\hat{\boldsymbol{\Sigma}}^{-1}(\boldsymbol{\mathcal{E}})\mathbf{Z})^{-1}\mathbf{D}'.$$

Then $\mathbf{M}(-\boldsymbol{\mathcal{E}}) = -\mathbf{M}(\boldsymbol{\mathcal{E}})$ and hence $\mathrm{E}(\mathbf{M}(\boldsymbol{\mathcal{E}})) = 0$. Thus, $\mathrm{E}(\hat{\boldsymbol{\gamma}}(\mathbf{Y})) = \mathbf{C}\boldsymbol{\Theta}\mathbf{D}'$. This completes the proof. □

## 3. Consistency

Since $\mathbf{Y}$ is associated with sample size $n$, we shall use $\mathbf{Y}_n$ to replace $\mathbf{Y}$ in (9) and then investigate the consistency of the estimator $\hat{\boldsymbol{\Sigma}}(\mathbf{Y}_n)$, as well as the consistency of the related estimator $\hat{\boldsymbol{\gamma}}(\mathbf{Y})$, as the sample size $n$ tends to infinity. Note that $\mathbf{X}$ and $\boldsymbol{\mathcal{E}}$ are also associated with the sample size $n$.

Recall that an estimator of $\boldsymbol{\Sigma}$ of the form $\mathbf{Y}_n' \mathbf{W}^* \mathbf{Y}_n$ is unbiased and invariant if and only if $\mathrm{tr}(\mathbf{W}^*) = 1$ and $\mathbf{W}^* \mathbf{X} = \mathbf{0}$; see Žežula [14]. Hence, the statistic $\hat{\boldsymbol{\Sigma}}(\mathbf{Y}_n) = \mathbf{Y}_n' \mathbf{W} \mathbf{Y}_n'$ defined in (6) is an unbiased and invariant estimator of $\boldsymbol{\Sigma}$ without the assumption of normality. Moreover, under the assumption of normality, the estimator $\hat{\boldsymbol{\Sigma}}(\mathbf{Y}_n)$ follows a Wishart distribution; see Hu [5].

Now, we shall investigate the consistency property of the estimator $\hat{\boldsymbol{\Sigma}}(\mathbf{Y}_n)$.

**Theorem 3.1.** *For model (1), the statistic $\hat{\boldsymbol{\Sigma}}(\mathbf{Y}_n)$ defined in (6) is a consistent estimator of the second-order parameter matrix $\boldsymbol{\Sigma}$.*

**Proof.** Since $\mathbf{Y}_n' \mathbf{W} \mathbf{Y}_n = (\mathbf{Y}_n - \mathbf{X}\boldsymbol{\Theta}\mathbf{Z}')' \mathbf{W} (\mathbf{Y}_n - \mathbf{X}\boldsymbol{\Theta}\mathbf{Z}')$, in the following discussion we can assume without loss of generality that $\mathbf{X}\boldsymbol{\Theta}\mathbf{Z}' = \mathbf{0}$. So, by (6),

$$\hat{\boldsymbol{\Sigma}}(\mathbf{Y}_n) = \frac{n}{n-m}\left(\frac{1}{n}\sum_{l=1}^n \boldsymbol{\mathcal{E}}_l \boldsymbol{\mathcal{E}}_l' - \frac{1}{n}\boldsymbol{\mathcal{E}}'\mathbf{P}_X \boldsymbol{\mathcal{E}}\right), \tag{11}$$

where $\boldsymbol{\mathcal{E}} = (\boldsymbol{\mathcal{E}}_1, \boldsymbol{\mathcal{E}}_2, \ldots, \boldsymbol{\mathcal{E}}_n)' \sim \mathcal{G}(\mathbf{0}, \mathbf{I}_n \otimes \boldsymbol{\Sigma})$.

Note that $(\boldsymbol{\mathcal{E}}_l \boldsymbol{\mathcal{E}}_l')_{l=1}^n$ is a random sample from a population with mean $\mathrm{E}(\boldsymbol{\mathcal{E}}_l \boldsymbol{\mathcal{E}}_l') = \boldsymbol{\Sigma}$. According to Kolmogorov's strong law of large numbers (see Rao [12], 2c.3 (iv)),

$$\frac{1}{n}\sum_{l=1}^n \boldsymbol{\mathcal{E}}_l \boldsymbol{\mathcal{E}}_l' \text{ converges almost surely to } \boldsymbol{\Sigma}. \tag{12}$$

Letting $\varepsilon > 0$, by Chebyshev's inequality and the fact that $\mathrm{E}(\mathbf{Y}'\mathbf{W}\mathbf{Y}) = \mathrm{tr}(\mathbf{W})\boldsymbol{\Sigma} + \mathrm{E}(\mathbf{Y})'\mathbf{W} \times \mathrm{E}(\mathbf{Y})$, we have

$$P\left(\left\|\frac{1}{\sqrt{n}}\mathbf{P}_X \boldsymbol{\mathcal{E}}\right\| \geq \varepsilon\right) \leq \frac{1}{n\varepsilon^2}\mathrm{E}[\mathrm{tr}(\boldsymbol{\mathcal{E}}'\mathbf{P}_X \boldsymbol{\mathcal{E}})] = \frac{1}{n\varepsilon^2}\mathrm{tr}(\mathrm{E}[\boldsymbol{\mathcal{E}}\boldsymbol{\mathcal{E}}']\mathbf{P}_X)$$

$$= \frac{1}{n\varepsilon^2}\mathrm{tr}(\mathbf{I}_n \mathrm{tr}(\boldsymbol{\Sigma})\mathbf{P}_X) = \frac{1}{n\varepsilon^2}\mathrm{tr}(\mathbf{P}_X)\mathrm{tr}(\boldsymbol{\Sigma}).$$

Since $\operatorname{tr}(\mathbf{P}_X) = \operatorname{rank}(\mathbf{X})$ is a constant, $P(\|\frac{1}{\sqrt{n}}\mathbf{P}_X\boldsymbol{\mathcal{E}}\| \geq \varepsilon)$ tends to $\mathbf{0}$ as the sample size $n$ tends to infinity. So,

$$\frac{1}{\sqrt{n}}\mathbf{P}_X\boldsymbol{\mathcal{E}} \text{ converges in probability to } \mathbf{0}. \tag{13}$$

Since convergence almost surely implies convergence in probability, by (12) and (13), we obtain from (11) that $\hat{\boldsymbol{\Sigma}}(\mathbf{Y}_n)$ converges to $\boldsymbol{\Sigma}$ in probability. This completes the proof. □

Now, we focus our attention on the consistency of the estimator $\hat{\boldsymbol{\gamma}}(\mathbf{Y}_n)$. We first prove the following lemma.

**Lemma 3.2.** $\mathbf{H}(\mathbf{Y}_n)$ *converges in probability to* $\mathbf{H}$*, where* $\mathbf{H}(\mathbf{Y}_n)$ *is defined in (8) and* $\mathbf{H} = \boldsymbol{\Sigma}^{-1}(\mathbf{P}_Z\boldsymbol{\Sigma}^{-1}\mathbf{P}_Z)^+$.

**Proof.** Note that the function $\mathbf{A}$ to $\mathbf{A}^+$ is not continuous. Since $\hat{\boldsymbol{\Sigma}}^{-1}(\mathbf{Y}_n)$ is positive definite with probability 1, by Lehmann [7], Lemma 5.3.2, and Theorem 3.1, we have

$$\hat{\boldsymbol{\Sigma}}^{-1}(\mathbf{Y}_n) \text{ converges in probability to } \boldsymbol{\Sigma}^{-1}. \tag{14}$$

Write

$$\mathbf{P}_Z = \mathbf{O}\boldsymbol{\Lambda}\mathbf{O}', \qquad \mathbf{Q}_n = \mathbf{O}'\hat{\boldsymbol{\Sigma}}(\mathbf{Y}_n)\mathbf{O},$$

where $\mathbf{O}$ is a $p \times p$ orthogonal matrix, $\boldsymbol{\Lambda} = \operatorname{diag}[\mathbf{0},\ \mathbf{I}_q]$ with $q = \operatorname{rank}(\mathbf{Z})$ and

$$\mathbf{Q}_n^{-1} = \mathbf{O}'\hat{\boldsymbol{\Sigma}}^{-1}(\mathbf{Y}_n)\mathbf{O} = \begin{bmatrix} \mathbf{G}_{11}(\mathbf{Y}_n) & \mathbf{G}_{12}(\mathbf{Y}_n) \\ \mathbf{G}_{21}(\mathbf{Y}_n) & \mathbf{G}_{22}(\mathbf{Y}_n) \end{bmatrix} = [\mathbf{G}_{ij}(\mathbf{Y}_n)]_{2\times 2}$$

with $\mathbf{G}_{22}(\mathbf{Y}_n)$ a $q \times q$ random matrix. By (14), for any $i,j = 1,2$, $\mathbf{G}_{ij}(\mathbf{Y}_n)$ converges in probability to $\mathbf{G}_{ij}$. Note that $(\mathbf{P}_Z\mathbf{C}\mathbf{P}_Z)^+ = \mathbf{P}_Z(\mathbf{P}_Z\mathbf{C}\mathbf{P}_Z)^+\mathbf{P}_Z$. Then

$$\begin{aligned}
\mathbf{H}(\mathbf{Y}_n) &= \hat{\boldsymbol{\Sigma}}^{-1}(\mathbf{Y}_n)(\mathbf{P}_Z\hat{\boldsymbol{\Sigma}}^{-1}(\mathbf{Y}_n)\mathbf{P}_Z)^+ = \mathbf{O}\mathbf{Q}_n^{-1}\boldsymbol{\Lambda}\mathbf{O}'(\mathbf{O}\boldsymbol{\Lambda}\mathbf{Q}_n^{-1}\boldsymbol{\Lambda}\mathbf{O}')^+\mathbf{O}\boldsymbol{\Lambda}\mathbf{O}' \\
&= \mathbf{O}(\mathbf{G}_{ij}(\mathbf{Y}_n))_{2\times 2}\boldsymbol{\Lambda}\mathbf{O}'(\mathbf{O}\boldsymbol{\Lambda}(\mathbf{G}_{ij}(\mathbf{Y}_n))_{2\times 2}\boldsymbol{\Lambda}\mathbf{O}')^+\mathbf{O}\boldsymbol{\Lambda}\mathbf{Q}' \\
&= \mathbf{O}\begin{bmatrix} \mathbf{0} & \mathbf{G}_{12}(\mathbf{Y}_n) \\ \mathbf{0} & \mathbf{G}_{22}(\mathbf{Y}_n) \end{bmatrix}\mathbf{O}'(\mathbf{O}\operatorname{diag}[\mathbf{0},\mathbf{G}_{22}(\mathbf{Y}_n)]\mathbf{O}')^+\mathbf{O}\boldsymbol{\Lambda}\mathbf{O}' \\
&= \mathbf{O}\begin{bmatrix} \mathbf{0} & \mathbf{G}_{12}(\mathbf{Y}_n) \\ \mathbf{0} & \mathbf{G}_{22}(\mathbf{Y}_n) \end{bmatrix}\operatorname{diag}[\mathbf{0},\mathbf{G}_{22}^{-1}(\mathbf{Y}_n)]\boldsymbol{\Lambda}\mathbf{O}' = \mathbf{O}\begin{bmatrix} \mathbf{0} & \mathbf{G}^*(\mathbf{Y}_n) \\ \mathbf{0} & \mathbf{I}_q \end{bmatrix}\mathbf{O}',
\end{aligned}$$

where $\mathbf{G}^*(\mathbf{Y}_n) = \mathbf{G}_{12}(\mathbf{Y}_n)\mathbf{G}_{22}^{-1}(\mathbf{Y}_n)$.

Similarly, $\mathbf{H}$ can be decomposed as

$$\mathbf{H} = \mathbf{O}\begin{bmatrix} \mathbf{0} & \mathbf{G}_{12}\mathbf{G}_{22}^{-1} \\ \mathbf{0} & \mathbf{I}_q \end{bmatrix}\mathbf{O}'.$$

Since $\mathbf{G}^*(\mathbf{Y}_n)$ converges in probability to $\mathbf{G}_{12}\mathbf{G}_{22}^{-1}$, we conclude that

$$\begin{bmatrix} \mathbf{0} & \mathbf{G}^*(\mathbf{Y}_n) \\ \mathbf{0} & \mathbf{I}_q \end{bmatrix} \text{ converges in probability to } \begin{bmatrix} \mathbf{0} & \mathbf{G}_{12}\mathbf{G}_{22}^{-1} \\ \mathbf{0} & \mathbf{I}_q \end{bmatrix},$$

namely, $\mathbf{H}(\mathbf{Y}_n)$ converges in probability to $\mathbf{H}$. This completes the proof. □

Based on Lemma 3.2, we obtain the following consistency result for the estimator $\hat{\boldsymbol{\gamma}}(\mathbf{Y}_n)$.

**Theorem 3.3.** *Assume that*

$$\lim_{n\to\infty} \frac{1}{n}(\mathbf{X}'\mathbf{X}) = \mathbf{R}, \tag{15}$$

*where $\mathbf{R}$ is a positive definite matrix. Then the statistic $\hat{\boldsymbol{\gamma}}(\mathbf{Y}_n)$ is a consistent estimator of the estimable parametric transformation $\boldsymbol{\gamma} = \mathbf{C}\boldsymbol{\Theta}\mathbf{D}'$.*

**Proof.** To prove that $\hat{\boldsymbol{\gamma}}(\mathbf{Y}_n)$ is a consistent estimator of $\boldsymbol{\gamma}$, by Slutsky's theorem (see Lehmann and Romano [8], Theorem 11.2.11), it suffices to show that $\hat{\boldsymbol{\Theta}}(\mathbf{Y}_n)$ is a consistent estimator of $\boldsymbol{\Theta}$.

Replacing $\mathbf{Y}$ with $\mathbf{X}\boldsymbol{\Theta}\mathbf{Z}' + \boldsymbol{\mathcal{E}}$ in (10), we decompose $\hat{\boldsymbol{\Theta}}(\mathbf{Y}_n)$ as $\mathbf{E}_n + \mathbf{F}_n$, where

$$\mathbf{E}_n = (\mathbf{X}'\mathbf{X})^{-1}\mathbf{X}'\mathbf{X}\boldsymbol{\Theta}\mathbf{Z}'\mathbf{H}(\mathbf{Y}_n)\mathbf{Z}(\mathbf{Z}'\mathbf{Z})^{-1}$$

and

$$\mathbf{F}_n = (\mathbf{X}'\mathbf{X})^{-1}\mathbf{X}'\boldsymbol{\mathcal{E}}\mathbf{H}(\mathbf{Y}_n)\mathbf{Z}(\mathbf{Z}'\mathbf{Z})^{-1}.$$

Note that $\mathbf{A}'\mathbf{A}(\mathbf{A}'\mathbf{A})^-\mathbf{B} = \mathbf{B}$ if $\mathbf{B}$ is estimable with respect to $\mathbf{A}$. Since $\mathbf{Z}'\mathbf{P}_Z = \mathbf{Z}'$, we obtain

$$\begin{aligned}\mathbf{E}_n &= \boldsymbol{\Theta}\mathbf{Z}'\hat{\boldsymbol{\Sigma}}^{-1}(\mathbf{Y}_n)(\mathbf{P}_Z\hat{\boldsymbol{\Sigma}}^{-1}(\mathbf{Y}_n)\mathbf{P}_Z)^+\mathbf{Z}(\mathbf{Z}'\mathbf{Z})^{-1}\\&= \boldsymbol{\Theta}\mathbf{Z}'\mathbf{P}_Z\hat{\boldsymbol{\Sigma}}^{-1}(\mathbf{Y}_n)\mathbf{P}_Z(\mathbf{P}_Z\hat{\boldsymbol{\Sigma}}^{-1}(\mathbf{Y}_n)\mathbf{P}_Z)^+\mathbf{Z}(\mathbf{Z}'\mathbf{Z})^{-1}\\&= \boldsymbol{\Theta}\mathbf{Z}'\mathbf{Z}(\mathbf{Z}'\mathbf{Z})^{-1} = \boldsymbol{\Theta}\end{aligned}$$

and

$$\mathbf{F}_n = n(\mathbf{X}'\mathbf{X})^{-1}\left(\frac{1}{\sqrt{n}}\mathbf{X}'\right)\left(\frac{1}{\sqrt{n}}\mathbf{P}_X\boldsymbol{\mathcal{E}}\right)\mathbf{H}(\mathbf{Y}_n)\mathbf{Z}(\mathbf{Z}'\mathbf{Z})^{-1}.$$

By (15), $\mathbf{X}'/\sqrt{n}$ are bounded. In fact, the elements of $\mathbf{X}'/\sqrt{n}$ are at most of order $n^{-1/2}$ (see the proof of Lemma 4.1 below). Then, by (13), (15) and Lemma 3.2, $\mathbf{F}_n$ converges in probability to $\mathbf{0}$. So, $\hat{\boldsymbol{\Theta}}(\mathbf{Y}_n)$ converges in probability to $\boldsymbol{\Theta}$. This completes the proof. □

## 4. Asymptotic normality

We investigated the consistency of the estimator $\hat{\boldsymbol{\gamma}}(\mathbf{Y}_n)$ in Section 3. In this section, we shall investigate the asymptotic normality of $\sqrt{n}[\hat{\boldsymbol{\gamma}}(\mathbf{Y}_n) - \boldsymbol{\gamma}]$.

First, we shall prove the following lemma.

**Lemma 4.1.** *Suppose that condition (15) holds. Let* $\mathbf{S} = (\mathbf{X}'\mathbf{X})^{-1}\mathbf{X}' = (\mathbf{s}_1, \mathbf{s}_2, \ldots, \mathbf{s}_n)_{m \times n}$, *where* $\mathbf{s}_l$ *is the* $l$*th column of* $(\mathbf{X}'\mathbf{X})^{-1}\mathbf{X}'$. *Then, for any* $l \in \{1, 2, \ldots, n\}$, *the* $m$ *elements of* $\sqrt{n}\mathbf{s}_l$ *are* $O(n^{-1/2})$.

**Proof.** Write $\mathbf{V} = \dfrac{1}{\sqrt{n}}\mathbf{X}' = [\mathbf{v}_1, \mathbf{v}_2, \ldots, \mathbf{v}_n]$. The transpose of $\mathbf{v}_l$ is an $m$-element row vector,

$$\mathbf{v}'_l = \left(\frac{1}{\sqrt{n}}x_{l1}, \frac{1}{\sqrt{n}}x_{l2}, \ldots, \frac{1}{\sqrt{n}}x_{lm}\right),$$

where $\mathbf{X} = [x_{ij}]_{n \times m}$. By (15), $\mathbf{V}\mathbf{V}' = n^{-1}\mathbf{X}'\mathbf{X}$ converges to a positive definite matrix $\mathbf{R}$. So, the elements of $\mathbf{V}\mathbf{V}' = \mathbf{v}_1\mathbf{v}'_1 + \mathbf{v}_2\mathbf{v}'_2 + \cdots + \mathbf{v}_n\mathbf{v}'_n$ are bounded. We claim that for any $l \in \{1, 2, \ldots, n\}$, the $m$ elements of $\mathbf{v}_l$ are all $O(n^{-1/2})$.

If this is not true, we can assume, without loss of generality, that one element of $\mathbf{v}_n$ is $O(n^{p-1/2})$ with $p > 0$. Then one element of $\mathbf{v}_n\mathbf{v}'_n$ would be $O(n^{2p-1})$. Hence, the corresponding element in matrix $\mathbf{V}\mathbf{V}' = \mathbf{v}_1\mathbf{v}'_1 + \mathbf{v}_2\mathbf{v}'_2 + \cdots + \mathbf{v}_n\mathbf{v}'_n$ would be $O(n^{2p})$, which is not bounded. This contradicts condition (15).

Note that

$$(\sqrt{n}\mathbf{s}_1, \sqrt{n}\mathbf{s}_2, \ldots, \sqrt{n}\mathbf{s}_n) = \sqrt{n}(\mathbf{X}'\mathbf{X})^{-1}\mathbf{X}' = n(\mathbf{X}'\mathbf{X})^{-1}\frac{1}{\sqrt{n}}\mathbf{X}'$$

$$= n(\mathbf{X}'\mathbf{X})^{-1}[\mathbf{v}_1, \mathbf{v}_2, \ldots, \mathbf{v}_n],$$

namely, for $l = 1, 2, \ldots, n$, $\sqrt{n}\mathbf{s}_l = n(\mathbf{X}'\mathbf{X})^{-1}\mathbf{v}_l$. Thus, for $l = 1, 2, \ldots, n$, the $m$ elements of $\sqrt{n}\mathbf{s}_l$ are also $O(n^{-1/2})$. This completes the proof. □

Now, we shall show the following important result on the asymptotic normality of $\sqrt{n}[\hat{\boldsymbol{\gamma}}(\mathbf{Y}_n) - \boldsymbol{\gamma}]$.

**Theorem 4.2.** *Under the assumption of condition (15),* $\sqrt{n}[\hat{\boldsymbol{\gamma}}(\mathbf{Y}_n) - \boldsymbol{\gamma}]$ *converges in distribution to* $\mathcal{N}_{s \times t}(\mathbf{0}, (\mathbf{C}\mathbf{R}^{-1}\mathbf{C}') \otimes (\mathbf{D}(\mathbf{Z}'\boldsymbol{\Sigma}^{-1}\mathbf{Z})^{-1}\mathbf{D}'))$.

**Proof.** First, by (8), we rewrite $\boldsymbol{\gamma}$ and $\hat{\boldsymbol{\gamma}}(\mathbf{Y}_n)$ as

$$\boldsymbol{\gamma} = \mathbf{C}\mathbf{S}\mathbf{X}\boldsymbol{\Theta}\mathbf{Z}'\mathbf{H}(\mathbf{Y}_n)\mathbf{K}\mathbf{D}'$$

and

$$\hat{\boldsymbol{\gamma}}(\mathbf{Y}_n) = \mathbf{C}\mathbf{S}\mathbf{Y}_n\mathbf{H}(\mathbf{Y}_n)\mathbf{K}\mathbf{D}',$$

where $\mathbf{K} = \mathbf{Z}(\mathbf{Z}'\mathbf{Z})^{-1}$. So,

$$\begin{aligned}
\hat{\boldsymbol{\gamma}}(\mathbf{Y}_n) - \boldsymbol{\gamma} &= \mathbf{CSY}_n\mathbf{H}(\mathbf{Y}_n)\mathbf{KD}' - \mathbf{CSX\Theta Z}'\mathbf{H}(\mathbf{Y}_n)\mathbf{KD}' \\
&= \mathbf{CS}(\mathbf{Y}_n - \mathbf{X\Theta Z}')\mathbf{H}(\mathbf{Y}_n)\mathbf{KD}' \\
&= \mathbf{CS}\boldsymbol{\mathcal{E}}\mathbf{H}(\mathbf{Y}_n)\mathbf{KD}' = \mathbf{CL}_n\mathbf{H}(\mathbf{Y}_n)\mathbf{KD}',
\end{aligned} \quad (16)$$

where $\mathbf{L}_n \equiv \mathbf{S}\boldsymbol{\mathcal{E}}$. Further, $\mathbf{L}_n$ is expressed as

$$\mathbf{L}_n = \sum_{l=1}^{n} \mathbf{s}_l \boldsymbol{\mathcal{E}}_l', \quad (17)$$

where $\mathbf{s}_l$ is the $l$th column vector of $\mathbf{S}$ and $\boldsymbol{\mathcal{E}}_l'$ is the $l$th row vector of the matrix $\boldsymbol{\mathcal{E}}$ with $\boldsymbol{\mathcal{E}} \sim \mathcal{G}(\mathbf{0}, \mathbf{I}_n \otimes \boldsymbol{\Sigma})$.

Next, we shall find the limiting distribution of $\sqrt{n}[\hat{\boldsymbol{\gamma}}(\mathbf{Y}_n) - \boldsymbol{\gamma}]$ through showing that

$$\sqrt{n}\mathbf{L}_n \text{ converges in distribution to } \mathcal{N}_{m \times p}(\mathbf{0}, \mathbf{R}^{-1} \otimes \boldsymbol{\Sigma}). \quad (18)$$

Since $\{\boldsymbol{\mathcal{E}}_l'\}_{l=1}^n$ are independent and identically distributed, for $\mathbf{t} \in \mathcal{M}_{p \times m}$, the characteristic function $\Psi_n(\mathbf{t})$ of $\sqrt{n}\mathbf{L}_n'$ is given by

$$\Psi_n(\mathbf{t}) = \mathrm{E}(\exp\{\mathrm{i}\,\mathrm{tr}(\mathbf{t}'\sqrt{n}\mathbf{L}_n')\}) = \mathrm{E}(\exp\{\mathrm{i}\,\mathrm{tr}(\mathbf{t}\sqrt{n}\mathbf{L}_n)\}) = \prod_{l=1}^{n} \Phi(\sqrt{n}\mathbf{t}\mathbf{s}_l),$$

where $\Phi(\cdot)$ is the characteristic function of $\boldsymbol{\mathcal{E}}_l'$.

Recall that for $u$ in the neighborhood of $0$,

$$\ln(1-u) = -u + f(u) \quad \text{with } f(u) = \tfrac{1}{2}u^2 + o(u^2). \quad (19)$$

If we write $p(u) = f(u)/u$, then from (19),

$$p(u) = o(u) \quad \text{as } u \to 0. \quad (20)$$

Also,

$$\Phi(\mathbf{x}) = 1 - \tfrac{1}{2}\mathbf{x}'\boldsymbol{\Sigma}\mathbf{x} + g(\mathbf{x}) \quad \text{for } \mathbf{x} \in \Re^m \quad \text{and} \quad g(\mathbf{x}) = o(\|\mathbf{x}\|^2) \quad \text{as } \mathbf{x} \to \mathbf{0}. \quad (21)$$

For $\varepsilon > 0$, there exists $\delta(\varepsilon) > 0$ such that

$$|g(\mathbf{x})| < \varepsilon\|\mathbf{x}\|^2 \quad \text{as } 0 < \|\mathbf{x}\| < \delta(\varepsilon). \quad (22)$$

Therefore, by (19) and (21), the characteristic function of $\sqrt{n}\mathbf{L}_n'$ can be decomposed as

$$\Psi_n(\mathbf{t}) = \exp\left\{\sum_{l=1}^{n} \ln(\Phi(\sqrt{n}\mathbf{t}\mathbf{s}_l))\right\}$$

$$= \exp\left\{\sum_{l=1}^{n} \ln\left(1 - \frac{n}{2}\mathbf{s}_l'\mathbf{t}'\mathbf{\Sigma}\mathbf{t}\mathbf{s}_l + g(\sqrt{n}\mathbf{t}\mathbf{s}_l)\right)\right\}$$
$$= \exp\left\{\sum_{l=1}^{n}\left[-\frac{1}{2}n\mathbf{s}_l'\mathbf{t}'\mathbf{\Sigma}\mathbf{t}\mathbf{s}_l + g(\sqrt{n}\mathbf{t}\mathbf{s}_l) + f\left(\frac{1}{2}n\mathbf{s}_l'\mathbf{t}'\mathbf{\Sigma}\mathbf{t}\mathbf{s}_l - g(\sqrt{n}\mathbf{t}\mathbf{s}_l)\right)\right]\right\} \quad (23)$$
$$= \exp\left\{-\frac{1}{2}\alpha_n + \beta_n + \eta_n\right\},$$

where

$$\alpha_n = \sum_{l=1}^{n} n\mathbf{s}_l'\mathbf{t}'\mathbf{\Sigma}\mathbf{t}\mathbf{s}_l = \operatorname{tr}\left(\sum_{l=1}^{n} n\mathbf{s}_l'\mathbf{t}'\mathbf{\Sigma}\mathbf{t}\mathbf{s}_l\right),$$
$$\beta_n = \sum_{l=1}^{n} g(\sqrt{n}\mathbf{t}\mathbf{s}_l)$$

and

$$\eta_n = \sum_{l=1}^{n} f(\tfrac{1}{2}n\mathbf{s}_l'\mathbf{t}'\mathbf{\Sigma}\mathbf{t}\mathbf{s}_l - g(\sqrt{n}\mathbf{t}\mathbf{s}_l)).$$

For $\alpha_n$, we have

$$\alpha_n = \operatorname{tr}(\mathbf{\Sigma}\mathbf{t}n\mathbf{S}\mathbf{S}'\mathbf{t}') = \operatorname{tr}(\mathbf{\Sigma}\mathbf{t}n(\mathbf{X}'\mathbf{X})^{-1}\mathbf{t}'). \quad (24)$$

By (15),

$$\lim_{n\to\infty}\alpha_n = \operatorname{tr}(\mathbf{R}^{-1}\mathbf{t}'\mathbf{\Sigma}\mathbf{t}) = (\operatorname{vec}(\mathbf{t}))'(\mathbf{R}^{-1}\otimes\mathbf{\Sigma})\operatorname{vec}(\mathbf{t}). \quad (25)$$

For $\beta_n$, by Lemma 4.1 and the continuity of $\mathbf{t}\mathbf{s}_l$, for the $\delta(\varepsilon) > 0$ in (22), there is an integer $N(\varepsilon) > 0$ such that for $n > N(\varepsilon)$,

$$0 < \|\sqrt{n}\mathbf{t}\mathbf{s}_l\| < \delta(\varepsilon) \quad \text{for all } l = 1,2,\ldots,n. \quad (26)$$

If we take $n > N(\varepsilon)$, then by (22) and (26),

$$|g(\sqrt{n}\mathbf{t}\mathbf{s}_l)| < \|\sqrt{n}\mathbf{t}\mathbf{s}_l\|^2\varepsilon. \quad (27)$$

So,

$$|\beta_n| < \sum_{l=1}^{n}\|\sqrt{n}\mathbf{t}\mathbf{s}_l\|^2\varepsilon = \varepsilon\operatorname{tr}(n\mathbf{t}\mathbf{S}\mathbf{S}'\mathbf{t}') = \varepsilon\operatorname{tr}(\mathbf{t}n(\mathbf{X}'\mathbf{X})^{-1}\mathbf{t}').$$

So, by (15), $\limsup_{n\to\infty}|\beta_n| \le \varepsilon\operatorname{tr}(\mathbf{t}\mathbf{R}^{-1}\mathbf{t}')$. Since $\varepsilon > 0$ is arbitrary, we obtain

$$\lim_{n\to\infty}\beta_n = 0. \quad (28)$$

For $\eta_n$, let
$$\lambda_l = \tfrac{1}{2}(\sqrt{n}\mathbf{t}\mathbf{s}_l)'\mathbf{\Sigma}(\sqrt{n}\mathbf{t}\mathbf{s}_l) - g(\sqrt{n}\mathbf{t}\mathbf{s}_l).$$

So, by (27),
$$|\lambda_l| < \tfrac{1}{2}(\sqrt{n}\mathbf{t}\mathbf{s}_l)'\mathbf{\Sigma}(\sqrt{n}\mathbf{t}\mathbf{s}_l) + \|\sqrt{n}\mathbf{t}\mathbf{s}_l\|^2\varepsilon. \tag{29}$$

Take $n > N(\varepsilon)$. By Lemma 4.1, the continuity of $\mathbf{t}\mathbf{s}_l$ and (20), increasing $N(\varepsilon)$ if necessary, we may suppose that for all $l$, $|p(\lambda_l)| < \varepsilon$. Since $f(\lambda_l) = p(\lambda_l)\lambda_l$,
$$|\eta_n| = \sum_{l=1}^{n}|f(\lambda_l)| = \sum_{l=1}^{n}|p(\lambda_l)||\lambda_l| \leq \sum_{l=1}^{n}\varepsilon|\lambda_l|.$$

So, by (29),
$$|\eta_n| \leq \sum_{l=1}^{n}\left[\frac{\varepsilon}{2}\operatorname{tr}(\sqrt{n}\mathbf{s}_l'\mathbf{t}'\mathbf{\Sigma}\mathbf{t}\sqrt{n}\mathbf{s}_l) + \|\sqrt{n}\mathbf{t}\mathbf{s}_l\|^2\varepsilon^2\right]$$
or
$$|\eta_n| \leq \sum_{l=1}^{n}\left[\frac{n\varepsilon}{2}\operatorname{tr}(\mathbf{t}'\mathbf{\Sigma}\mathbf{t}\mathbf{s}_l\mathbf{s}_l') + \operatorname{tr}(\sqrt{n}\mathbf{t}\mathbf{s}_l(\sqrt{n}\mathbf{t}\mathbf{s}_l)'\varepsilon^2)\right],$$
namely,
$$|\eta_n| \leq \frac{\varepsilon}{2}\operatorname{tr}(\mathbf{t}'\mathbf{\Sigma}\mathbf{t}n\mathbf{S}\mathbf{S}') + \operatorname{tr}(\mathbf{t}n\mathbf{S}\mathbf{S}'\mathbf{t}')\varepsilon^2. \tag{30}$$

Note that $n\mathbf{S}\mathbf{S}' = n(\mathbf{X}'\mathbf{X})^{-1}$. Since $\varepsilon$ is arbitrary, by (15) and (30),
$$\lim_{n\to\infty}\eta_n = 0. \tag{31}$$

By (25), (28) and (31), we obtain from (23) that
$$\lim_{n\to\infty}\Psi_n(\mathbf{t}) = \exp\{-\tfrac{1}{2}(\operatorname{vec}(\mathbf{t}))'(\mathbf{R}^{-1}\otimes\mathbf{\Sigma})\operatorname{vec}(\mathbf{t})\}. \tag{32}$$

So, by Lévy's continuity theorem, $\sqrt{n}\mathbf{L}_n$ in (17) converges in distribution to the normal distribution $\mathcal{N}_{m\times p}(\mathbf{0},\mathbf{R}^{-1}\otimes\mathbf{\Sigma})$, as was claimed in (18).

Finally, by Lemma 3.2, (16), (18) and Muirhead [9], Theorem 1.2.6, we obtain that

$\sqrt{n}[\hat{\boldsymbol{\gamma}}(\mathbf{Y}_n) - \boldsymbol{\gamma}]$ converges in distribution to $\mathcal{N}_{s\times t}(\mathbf{0},\mathbf{C}\mathbf{R}^{-1}\mathbf{C}'\otimes(\mathbf{D}\mathbf{K}'\mathbf{H}'\mathbf{\Sigma}\mathbf{H}\mathbf{K}\mathbf{D}'))$.

Replacing $\mathbf{H}$ and $\mathbf{K}$ with $\mathbf{\Sigma}^{-1}(\mathbf{P}_Z\mathbf{\Sigma}^{-1}\mathbf{P}_Z)^+$ and $\mathbf{Z}(\mathbf{Z}'\mathbf{Z})^{-1}$, respectively, we conclude that

$\sqrt{n}[\hat{\boldsymbol{\gamma}}(\mathbf{Y}_n) - \boldsymbol{\gamma}]$ converges in distribution to $\mathcal{N}_{s\times t}(\mathbf{0},(\mathbf{C}\mathbf{R}^{-1}\mathbf{C}')\otimes(\mathbf{D}\mathbf{T}\mathbf{D}'))$,

where $\mathbf{T} = (\mathbf{Z}'\mathbf{Z})^{-1}\mathbf{Z}'(\mathbf{P}_Z\mathbf{\Sigma}^{-1}\mathbf{P}_Z)^+\mathbf{Z}(\mathbf{Z}'\mathbf{Z})^{-1} = (\mathbf{Z}'\mathbf{\Sigma}^{-1}\mathbf{Z})^{-1}$ (see (3)). Thus, the proof is complete. □

Under model (1), hypotheses of the form

$$\mathrm{H}: \boldsymbol{\gamma} \equiv \mathbf{C}\boldsymbol{\Theta}\mathbf{D}' = \mathbf{0}$$

are usually considered; see Potthoff and Roy [11].

From Theorem 4.2 and Slutsky's theorem, the following corollary provides the asymptotic behavior of $\sqrt{n}\hat{\boldsymbol{\gamma}}(\mathbf{Y}_n)$ under H.

**Corollary 4.3.** *Under the assumption of condition (15), if matrices $\mathbf{C}(\mathbf{X}'\mathbf{X})^{-1}\mathbf{C}'$ and $\mathbf{D}(\mathbf{Z}'\hat{\boldsymbol{\Sigma}}^{-1}(\mathbf{Y}_n)\,\mathbf{Z})^{-1}\mathbf{D}'$ are non-singular, then the statistic $(\mathbf{C}n(\mathbf{X}'\mathbf{X})^{-1}\mathbf{C}')^{-1/2}\sqrt{n}\hat{\boldsymbol{\gamma}}(\mathbf{Y}_n)(\mathbf{D}(\mathbf{Z}' \times \hat{\boldsymbol{\Sigma}}^{-1}(\mathbf{Y}_n)\mathbf{Z})^{-1}\mathbf{D}')^{-1/2}$ under H converges in distribution to $\mathcal{N}_{s\times t}(\mathbf{0}, \mathbf{I})$.*

*Remark 4.4.* Lemma 2.1 tells us that $\hat{\boldsymbol{\gamma}}(\mathbf{Y})$ is an unbiased estimator of $\boldsymbol{\gamma}$ under the assumption of $\mathcal{E}$ being symmetric about the origin. In general, it is very difficult to obtain the covariance matrix of $\hat{\boldsymbol{\gamma}}(\mathbf{Y})$, even under the assumption of normality. However, under condition (15), Theorem 4.2 gives us an approximate covariance matrix $(\mathbf{C}(\mathbf{X}'\mathbf{X})^{-1}\mathbf{C}') \otimes (\mathbf{D}(\mathbf{Z}'\hat{\boldsymbol{\Sigma}}^{-1}(\mathbf{Y})\mathbf{Z})^{-1}\mathbf{D}')$ of $\hat{\boldsymbol{\gamma}}(\mathbf{Y})$ for large sample size $n$, without the assumption of normality.

We now conclude this paper by discussing the example in Potthoff and Roy [11]. No assumption of normality is made in our discussion.

*Example 4.5.* There are $m$ groups of animals, with $r$ animals in the $j$th group and each group being subjected to a different treatment. Animals in all groups are measured at the same $p$ time points, $t_1, t_2, \ldots, t_p$. The observations of different animals are independent, but the $p$ observations on each animal are assumed to have a covariance matrix $\boldsymbol{\Sigma}$.

Based on the problem and our discussion, $m$ remains constant, while $r$ tends to infinity. For $i = 1, 2, \ldots, m$, the growth curve associated with the $i$th group is

$$\theta_{i0} + \theta_{i1}x + \theta_{i2}x^2 + \cdots + \theta_{1q-1}x^{q-1}.$$

Put

$$\mathbf{X} = (\mathbf{x}_1, \mathbf{x}_2, \ldots, \mathbf{x}_m),$$

where $\mathbf{x}_i = [\delta_{1i}\mathbf{e}'_r, \delta_{2i}\mathbf{e}'_r, \ldots, \delta_{ti}\mathbf{e}'_r]'$, $\mathbf{e}_r = (1, 1, \ldots, 1)' \in \Re^r$, $\delta_{ij}$ are the Kronecker symbols, $n = rm$,

$$\boldsymbol{\theta}_i = (\theta_{i0}, \theta_{i1}, \theta_{i2}, \ldots, \theta_{i\,q-1}), \qquad \boldsymbol{\Theta} = (\boldsymbol{\theta}'_1, \boldsymbol{\theta}'_2, \ldots, \boldsymbol{\theta}'_m)'$$

and

$$\mathbf{Z} = \begin{bmatrix} 1 & t_1 & t_1^2 & \cdots & t_1^{q-1} \\ 1 & t_2 & t_2^2 & \cdots & t_2^{q-1} \\ . & . & . & \cdots & . \\ 1 & t_p & t_p^2 & \cdots & t_p^{q-1} \end{bmatrix}.$$

The observation data matrix $\mathbf{Y}_n$ can be written as

$$\mathbf{Y}_n = \mathbf{X}\boldsymbol{\Theta}\mathbf{Z}' + \boldsymbol{\mathcal{E}},$$

where $\boldsymbol{\mathcal{E}} = (\boldsymbol{\mathcal{E}}_1, \boldsymbol{\mathcal{E}}_2, \ldots, \boldsymbol{\mathcal{E}}_n)'$ with $\boldsymbol{\mathcal{E}}_1, \boldsymbol{\mathcal{E}}_2, \ldots, \boldsymbol{\mathcal{E}}_n$ being independent and identically distributed with mean $\mathbf{0}$ and covariance $\boldsymbol{\Sigma}$. Then, by (14),

$$\mathbf{R}^{-1} = \lim_{r \to \infty} n(\mathbf{X}'\mathbf{X})^{-1} = m\mathbf{I}.$$

By (10),

$$\hat{\boldsymbol{\Theta}}(\mathbf{Y}_n) = \frac{m}{n}\mathbf{X}'\mathbf{Y}_n\hat{\boldsymbol{\Sigma}}^{-1}(\mathbf{Y}_n)(\mathbf{P}_Z(\hat{\boldsymbol{\Sigma}}^{-1}(\mathbf{Y}_n))\mathbf{P}_Z)^+\mathbf{K}.$$

For the estimable parametric transformation of the form $\boldsymbol{\gamma} = \mathbf{C}\boldsymbol{\Theta}\mathbf{D}'$ with given $\mathbf{C} \in \mathscr{M}_{s \times m}$ and $\mathbf{D} \in \mathscr{M}_{t \times q}$, the two-stage generalized least-squares estimator is given by

$$\hat{\boldsymbol{\gamma}}(\mathbf{Y}_n) = \mathbf{C}\hat{\boldsymbol{\Theta}}(\mathbf{Y}_n)\mathbf{D}'.$$

It follows from Theorem 4.2 that $\sqrt{n}[\hat{\boldsymbol{\gamma}}(\mathbf{Y}_n) - \boldsymbol{\gamma}]$ converges in distribution to the normal distribution $\mathcal{N}_{s \times t}(\mathbf{0}, (m\mathbf{C}\mathbf{C}') \otimes (\mathbf{D}(\mathbf{Z}'\boldsymbol{\Sigma}^{-1}\mathbf{Z})^{-1}\mathbf{D}'))$.

Moreover, if we try to test that all $m$ growth curves are equal, except possibly for the additive constant $\theta_{i0}$, then we take $\mathbf{C}$ to be a matrix whose last column contains all $-1$'s and whose first $(m-1)$ columns constitute the identity matrix, and $\mathbf{D}$ to be a $(q-1) \times q$ matrix whose first column contains all 0's and whose last $(q-1)$ columns constitute the identity matrix, namely, taking

$$\mathbf{C} = [\mathbf{I}_{m-1} \quad -\mathbf{1}_{m-1}]_{(m-1) \times m}, \qquad \mathbf{D} = [\mathbf{0} \quad \mathbf{I}_{q-1}]_{(q-1) \times q},$$

where $\mathbf{1}_{m-1} = (1, 1, \ldots, 1)'$, and hypothesis $\mathrm{H}_0 : \mathbf{C}\boldsymbol{\Theta}\mathbf{D}' = \mathbf{0}$. Obviously, matrices $\mathbf{C}(\mathbf{X}'\mathbf{X})^{-1}\mathbf{C}'$ and $\mathbf{D}(\mathbf{Z}'\hat{\boldsymbol{\Sigma}}^{-1}(\mathbf{Y}_n)\mathbf{Z})^{-1}\mathbf{D}'$ are non-singular. It follows from Corollary 4.3 that statistic $(\mathbf{C}n(\mathbf{X}'\mathbf{X})^{-1}\mathbf{C}')^{-1/2}\sqrt{n}\hat{\boldsymbol{\gamma}}(\mathbf{Y}_n)(\mathbf{D}(\mathbf{Z}'\hat{\boldsymbol{\Sigma}}^{-1}(\mathbf{Y}_n)\mathbf{Z})^{-1}\mathbf{D}')^{-1/2}$ under $\mathrm{H}_0$ converges in distribution to $\mathcal{N}_{(m-1) \times (q-1)}(\mathbf{0}, \mathbf{I})$.

## Acknowledgements

The authors would like to thank Professor C.S. Wong for his guidance and encouragement in obtaining these results. The authors are also deeply grateful to the anonymous referee and the Associate Editor for their helpful suggestions which led to the improved version of this paper.